\documentclass{amse-new}
\numberwithin{equation}{section} %%% Equations numbered by section. If you don't want it, please delete it.

\begin{document}

 \PageNum{1}
 \Volume{201x}{Sep.}{x}{x}
 \OnlineTime{August 15, 201x}
 \DOI{0000000000000000}
 \EditorNote{Received x x, 201x, accepted x x, 201x}

\abovedisplayskip 6pt plus 2pt minus 2pt \belowdisplayskip 6pt
plus 2pt minus 2pt
%%%%%%%%%%%%%%%%
\def\vsp{\vspace{1mm}}
\def\th#1{\vspace{1mm}\noindent{\bf #1}\quad}
\def\proof{\vspace{1mm}\noindent{\it Proof}\quad}
\def\no{\nonumber}
\newenvironment{prof}[1][Proof]{\noindent\textit{#1}\quad }
{\hfill $\Box$\vspace{0.7mm}}
\def\q{\quad} \def\qq{\qquad}
\allowdisplaybreaks[4]
%%%%%%%%%%%%%%%%%%%%%%%%%%%%%%%%%%%%%%%%%%%%%%%%%%%%%%%%%%%%%%%%%%%%%%%%%%%%%%%%%%%%%%%%%%%%%%%
%%-------------------       Beginning of  Author's Definitions       -------------------%%
%%                     Note: You may add your own definitions here.

%%-------------------         the end of  Author's Definitions           -------------------%%

\AuthorMark{Zhang Y.S.}                             %%%  appear on the head of even pages  %%%

\TitleMark{On Lawson's area-minimizing hypercones}  %%% Running Title, appear on the head of odd pages  %%%

\title{On Lawson's area-minimizing hypercones      %%%   Main Title of your paper  %%%
\footnote{
Partially sponsored by the Fundamental Research Funds for the Central Universities, the SRF for ROCS, SEM{,}
%%%%%%%%%%%%%%%%%%%%%%%%%%%%%%%%%%%%%%%%%%%%%%%%%%%%%%%%%%%%%%%%%%%%%%   20160922
%\textcolor{red}
{
NSFC (Grant Nos. 11526048, 11601071)
}
%%%%%%%%%%%%%%%%%%%%%%%%%%%%%%%%%%%%%%%%%%%%%%%%%%%%%%%%%%%%%%%%%%%%%%    20160922
and 
the NSF (Grant No. 0932078 000)
while the author was in residence at the MSRI during the 2013 Fall.}}            %%%   the Fund which you are supported by  %%%

\author{Yongsheng \uppercase{Zhang}}             %%%  1st Author's information   %%%
    {School of Mathematics and Statistics, Northeast Normal University\\ ChangChun 130024, P.R. China
    \\E-mail\,$:$ yongsheng.chang@gmail.com}

\maketitle%

\Abstract{We show the area-minimality property of all homogeneous area-minimizing hypercones in Euclidean spaces (classified by Lawlor)
following Lawson's original idea in his 72' Trans. A.M.S. paper ``The equivariant Plateau problem and interior regularity".
Moreover, each of them enjoys (coflat) calibrations singular only at the origin.}      % the abstract

        % the keywords

\MRSubClass{53C38, 28A75}      % MR(2000) Subject Classification

  \section{Introduction}\label{Section1} 
                    Let $X$ be a Riemannian manifold and $G$ a compact, connected group of isometries of $X$.
                    In \cite{WHL} Hsiang and Lawson showed that
                       a $G$-invariant submanifold $N$ of cohomogeneity $k$ in $X$ is minimal if and only if $N^*/G$ is minimal in $X^*/G$
                       with respect to a natural metric on the orbit space.
                   Here $*$ means the union of principal orbits under action $G$.
                   In particular they classified compact homogeneous minimal hypersurfaces in standard spheres. 
                   
                   In this paper $X=\mathbb R^{n+1}$ with Euclidean metric and $M$ 
                   is a submanifold in the unit sphere $S^{n}$.
                   We call
                            $$
                            C(M)=\{tx: 0\leq t <+\infty, x\in M\}
                            $$
                  the cone over $M$ and $M$ the link of $C(M)$.
                  It follows that $M$ is minimal in $S^{n}$ if and only if $C(M)\sim0$ is minimal in $\mathbb R^{n+1}$.
                  
                  In \cite{BL} Lawson first proved
                           that there {\it always} exists a $G$-invariant solution to the Plateau problem for any $G$-invariant boundary $M$ of codimension 2, 
                  and then he considered further which homogeneous minimal hypersurfaces in their classification give rise to area-minimizing cones,
                  in the sense of Federer, namely the truncated cone by the unit ball centered at the origin is an area-minimizer among all integral currents 
{(introduced in \cite{FF})}  %%%%%%%%%%%%%%%%%%%%%%%%%%%%%%%%%%%%%%%%%%        20160922
                  sharing the boundary $M$ (these hypercones are closely related to the Bernstein problem, see \cite{BDG, LS}).
                  The question turns out to be equivalent to whether the line segment between the ``0"-point and ``$M$"-point in the orbit space is minimizing.
                  He gained the following variety of minimizing hypercones by means of constructing calibration $1$-forms on orbit spaces.

      \begin{theorem}[Lawson \cite{BL}]\label{th:1.1}
       There exist link manifolds $M$ of types
       $ \mathbf{(I)}$ $Sp(3)/Sp(1)^3$ in $\mathbb R^{14}$,
        $F_4/Spin(8)$ in $\mathbb R^{26}$,
         and moreover,
        $\mathbf{(II)}$ $S^{r-1}\times S^{s-1}$ in $\mathbb R^{r+s}$ for $r,s\geq 2$
        with $r+s\geq 10$ or with $r+s=9$ and $|r-s|\leq 5$ or with $r=s=4$,
        $SO(2)\times SO(k)/\mathbb Z_2\times SO(k-2)$ in $\mathbb R^{2k}$ for $k\geq 10$,
        $SU(2)\times SU(k)/T^1\times S
%\textcolor{red}
{U}      %%%%%%%%%%%%%%%%%%%%%%%%%%%%%%%%%%%%%%%%%%        20160922
        (k-2)$ in $\mathbb R^{4k}$ for $k\geq 5$,
        %$SU(5)/SU(2)\times SU(2)$ in $\mathbb R^{20}$,
%\textcolor{red}
{$U(5)/SU(2)\times SU(2) \times T^1$} in $\mathbb R^{20}$,           %%%%%%%%%%%%%%%%%%%%%%%%%%%%%%%%%%%%%%%%%%        20160922
        $U(1)\cdot Spin(10)/T^1\cdot SU(4)$ in $\mathbb R^{32}$,
        and 
        $Sp(2)\times Sp(k)/Sp(1)^2\times Sp(k-2)$ in $\mathbb R^{8k}$ for $k\geq 2$,
        such that $C(M)$ is area-minimizing.        
      \end{theorem}
       
       \begin{remark}
       Here we put overlapping types of  $ \mathbf{(I)}$ and $ \mathbf{(II)}$ in \cite{BL} into class $ \mathbf{(II)}$.
       \end{remark}

       \begin{remark}
       For example, cones $C_{m,n}\triangleq C(S^{m}(\sqrt{\frac{m}{m+n}})\times S^n(\sqrt{\frac{n}{m+n}}))$
       where $m=r-1$ and $n=s-1$ for $r,s$ in the theorem are area-minimizing.
       \end{remark}
       
        From the point of view of calibrated geometry, the calibration on the orbit space
        induces a coflat calibration in the Euclidean space (see \S\ref{s4})
        that is singular in a set of codimension at least two.
Hence,
%\textcolor{red}{,}          %%%%%%%%%%%%%%%%%%%%%%%%%%%%%%%%%%%%%%%%%%        20160922
        by the coflat version of fundamental theorem of calibrated geometry (Theorem 4.9 in \cite{HL2}), homogeneous hypercones with such calibrations
        are area-minimizing.
        In fact similar idea also can work for cones of higher codimension,
        for example, see \cite{Ch}.

       By studying Lawson's calibrations,
       we give a proof of the {\it affirmative} part of the table on page 86
       (classification of homogeneous area-minimizing hypercones) in \cite{Law} by verifying the rest cases.
       
       \begin{theorem}\label{th:1.4}
       $M$ of class $\mathbf{(II)}$ in the above theorem can also be of type $S^{2}\times S^{4}$ in $\mathbb R^{8}$,
         $SO(2)\times SO(k)/\mathbb Z_2\times SO(k-2)$ in $\mathbb R^{2k}$ for $k=9$,
and $SU(2)\times SU(k)/T^1\times S{U}(k-2)$ in $\mathbb R^{4k}$ for $k=4$,  %%%%%%%%%%%%%%%%%%%%%%%%%%%%        20160922
        such that $C(M)$ is area-minimizing.
        Moreover, types of $S^{1}\times S^{5}$ and $SO(2)\times SO(8)/\mathbb Z_2\times SO(6)$
        correspond to stable minimal hypercones.     
      \end{theorem}
         
         \begin{remark}
         The area-minimality of $C_{2,4}$ in $\mathbb R^8$ was given by Simoes \cite{PS1,PS2}.
         \end{remark}
         
         According to Hardt and Simon \cite{HS}, here comes a corollary.
         
         \begin{corollary}
         All homogeneous area-minimizing hypercones
          are strictly area-minimizing.  
       \end{corollary}
 
       \begin{remark}
       The strict area-minimality of $C_{2,4}$ was first proved by Lin \cite{FL}.
       \end{remark}
       
        We obtain the next theorem
       through a particular kind of modifications on Lawson's calibrations.

       \begin{theorem}\label{th:1.8}
       Each area-minimizing $C_{r-1,s-1} \subset \mathbb R^{r+s}$ has a coflat calibration singular only at the origin.
       \end{theorem}  
       
       An observation on
       the calibration O.D.E. \eqref{Code} leads to a more general result.
       
       \begin{theorem}\label{th:1.9}
       Every homogeneous area-minimizing hypercone has a coflat calibration singular only at the origin.
       \end{theorem}
       
       It may be interesting to consider whether every area-minimizing hypercone has such nice calibrations.
      See \cite{Z12} for applications of these calibrations. % \\{\ }     
%      \\
%      {\ }
  \\

  %%%%%%%%%%%%%%%%%%%%%%%%%%%%%%%%%%%%%

        \section{Preliminaries}     

         Suppose a {connected} compact  Lie group $G$ has an orthogonal representation on $X$ of cohomogeneity 2
         and $M_0$ is the principal orbit of greatest volume in $S^n$. 
         Then a celebrated classification due to \cite{WHL} and \cite{TT} is the following
{table,} 
         based on which much interesting research has been done, for example, recent papers \cite{TXY} and \cite{MO}.
          Here we only list the needed information and refer readers to their original papers for details.
% \\{\ }\\
\begin{table}
\begin{center}
{
\begin{tabular}{r|c|c|c|c}   
        &   $G$                                 &                  type of $M_0$                                                                        &                        $A$             &                   $V^2$ \\[0.5ex]\hline\hline
1     &    $SO(r)\times SO(s)$     &             $S^{r-1}\times S^{s-1}\ \ r,s\geq2$                                                             &                        $\frac{\pi}{2}$                &                    $c\cdot x^{2r-2}y^{2s-2}$\\\hline
2     &    $SO(2)\times SO(k)$   &             $\frac{SO(2)\times SO(k)}{\mathbb Z_2\times SO(k-2)}\ \ k\geq3$   &                         $\frac{\pi}{4}$               &                    $c\cdot (xy)^{2k-4}(x^2-y^2)^2$\\\hline
3     &    $SU(2)\times SU(k)$   &             $\frac{SU(2)\times SU(k)}{T^1\times SU(k-2)}\ \ k\geq2$                    &                         $\frac{\pi}{4}$               &                    $c\cdot (xy)^{4k-6}(x^2-y^2)^4$\\\hline
4     &    $Sp(2)\times Sp(k)$    &             $\frac{Sp(2)\times Sp(k)}{Sp(1)^2\times Sp(k-2)}\ \ k\geq2$               &                         $\frac{\pi}{4}$               &                    $c\cdot (xy)^{8k-10}(x^2-y^2)^8$\\\hline
5     &                           $U(5)$    &             $\frac{U(5)}{SU(2)\times SU(2)\times T^1}$                           &                         $\frac{\pi}{4}$               &                    $c\cdot (xy)^{2}$Im$\{(x+iy)^4\}^8$\\\hline
6  & $U(1)\cdot Spin(10)$    &             $\frac{U(1)\cdot Spin(10)}{T^1\cdot SU(4)}$                       &                         $\frac{\pi}{4}$               &                    $c\ (xy)^6$Im$\{(x+iy)^4\}^{12}$
\\\hline
7    &                        $SO(3)$    &             $\frac{SO(3)}{\mathbb Z_2+ \mathbb Z_2}$                                   &                         $\frac{\pi}{3}$               &                    $c\ \cdot\ $Im$\{(x+iy)^3\}^2$\\\hline
8     &                        $SU(3)$    &             $\frac{SU(3)}{T^2}$                                                                      &                         $\frac{\pi}{3}$               &                    $c\ \cdot\ $Im$\{(x+iy)^3\}^4$\\\hline
9     &                        $Sp(3)$     &             $\frac{Sp(3)}{Sp(1)^3}$                                                               &                         $\frac{\pi}{3}$               &                    $c\ \cdot\ $Im$\{(x+iy)^3\}^8$\\\hline
10  &                        $F_4$         &             $\frac{F_4}{Spin(8)}$                                                                    &                         $\frac{\pi}{3}$               &                    $c\ \cdot\ $Im$\{(x+iy)^3\}^{16}$\\\hline
11    &                        $Sp(2)$      &             $\frac{Sp(2)}{T^2}$                                                                       &                         $\frac{\pi}{4}$               &                    $c\ \cdot\ $Im$\{(x+iy)^4\}^4$\\\hline
12  &                        $G_2$        &             $\frac{G_2}{T^2}$                                                                          &                         $\frac{\pi}{6}$               &                    $c\ \cdot\ $Im$\{(x+iy)^6\}^4$\\\hline

13  & $SO(4)$    &             $\frac{SO(4)}{\mathbb Z_2+ \mathbb Z_2}$                       &                         $\frac{\pi}{6}$               &                    $c\ \cdot\ $Im$\{(x+iy)^6\}^2$\\\hline
\hline
\end{tabular}
}
\caption{Minimal homogeneous hypersurfaces in spheres}
\end{center}
\end{table}
%{\ }

Each orbit space is a cone of angle $A$ in $\mathbb R^2$ endowed with metric $$ds^2=V^2(dx^2+dy^2),$$
where $V$ is the volume function of orbits. %(We ignore the constant factor $c$ of $V^2$.) 
Let $\pi$ be the projection map to orbit space. Then it can be shown that (up to a constant factor) the volume of any $G$-invariant hypersurface $N$ equals the length of $\pi (N)$.

In order to search for area-minimizing cone $C(M_0)$,
Lawson made the following try for class \textbf{(II)}, i.e., Row $1-6$ in the table.

         First note that $A=\frac{\pi}{4}$ for Row $2-6$.
         Consider $re^{i\theta}=z^2$ with $z=x+yi$.
         Now the orbit space is stretched to a cone $Q$ of angle $\frac{\pi}{2}$ with metrics %(up to a scalar factor) 
                $r^{2k-3}s^{2k-4} c^2 (r^2d\theta^2+dr^2)$,
                 $r^{4k-3}s^{4k-6}c^4 (r^2d\theta^2+dr^2)$,
                  $r^{8k-3}s^{8k-10} c^8 (r^2d\theta^2+dr^2)$,
                   $r^{17}s^{10} c^8 (r^2d\theta^2+dr^2)$
                   and
                   $r^{29}s^{18} c^{12} (r^2d\theta^2+dr^2)$ respectively.
          Here $c=\cos\theta$ and $s=\sin \theta$.
          For metric $ds^2=r^l c^ps^q(r^2d\theta+dr^2)$ on $Q$ (including Row 1),
          the angle $\theta_0$ corresponding to $M_0$ satisfies $\tan^2 \theta_0=q/p$.
          Set $\tau=\cos^p\theta_0\sin^q\theta_0.$
          Then, under $d\tilde s^2=\frac{1}{\tau}ds^2$ on $Q$, Lawson considered the function of type
             \begin{equation}\label{0}
              f =\alpha^{-1}r^\alpha \phi^\beta,
              \end{equation}
where
$2\alpha=l+2,\ \phi=\frac{c^ps^q}{\tau}\  \mathrm{and}\ \beta \mathrm{\ is\ some\ real\ number} \mathrm{\ to\ be\ determined}.$
Since 
\[
\phi'(\theta)=\phi(\theta)(q\cot \theta-p\tan \theta),
\]
it follows
$$df=r^{\alpha-1}\phi^\beta dr+(\frac{\beta}{\alpha})r^\alpha\phi^\beta(q\cot\theta-p\tan\theta)d\theta.$$

Note that for a $1$-form its comass (see \cite{HL2}) and usual norm coincide.
Therefore if for some $\beta$
\begin{equation}\label{code}
r^{2\alpha-2}\phi^{2\beta} [1+(\frac{\beta}{\alpha})^2(q\cot\theta-p\tan\theta)^2]\leq \frac{r^lc^ps^q}{\tau},
\end{equation}
or simply
\begin{equation}\label{1}
\phi^{2\beta-1} [1+(\frac{\beta}{\alpha})^2(q\cot\theta-p\tan\theta)^2]\leq1,
\end{equation}
then $df$ is a calibration $1$-form, smooth on $Q^\circ$ and continuous to $\partial Q$, calibrating the ray $\theta=\theta_0$,
and hence $C(M_0)$ is area-minimizing.
Denote the left-hand side of \eqref{1} by $\psi$.
It is clear that $\psi$ is pointwise the square of comass of $df$ with respect to $d\tilde s^2$.

               %%%%%%%%%%%%%%%%%%%
           {\ }\\
             
               \section{Proof of Theorem \ref{th:1.4}}
Consider the derivative of $\psi$.
                \[
                \begin{split}
                 \psi'                 =\phi^{2\beta-1}(q\cot\theta-p\tan\theta)          \Bigg\{&(2\beta-1)\left[1+(\frac{\beta}{\alpha})^2(q\cot\theta-p\tan\theta)^2\right]\\
                                                                                                           &+2(\frac{\beta}{\alpha})^2(q\cot\theta-p\tan\theta)'\Bigg\}\\
                 \end{split}
                 \]
{\ }\\          %%%%%%%%%%%%%%%%%%%%%%%%%%%%%%%%%%%%%%%%%%        20160922         
{L}et $\eta_\beta$ stand for the part in braces. Then
                  \[
                \begin{split}
                 \eta_\beta                =&(2\beta-1)\left[1+(\frac{\beta}{\alpha})^2(q\cot\theta-p\tan\theta)^2\right]
                                                                                                           +2(\frac{\beta}{\alpha})^2(-\frac{q}{s^2}-\frac{p}{c^2})\\
                                        =&(2\beta-1)\left[1+(\frac{\beta}{\alpha})^2(q^2\cot^2\theta+p^2\tan^2\theta-2pq)\right]+2(\frac{\beta}{\alpha})^2(-\frac{q}{s^2}-\frac{p}{c^2})\\
                                        =&(2\beta-1)-2(\frac{\beta}{\alpha})^2\left[(2\beta-1)pq+p+q\right]\\
                                        &+(\frac{\beta}{\alpha})^2
                                        \Big\{\left[(2\beta-1)q^2-2q\right]\cot^2\theta+\left[(2\beta-1)p^2-2p\right]\tan^2\theta\Big\}
                 \end{split}
                 \]
{\ } \\     %%%%%%%%%%%%%%%%%%%%%%%%%%%%%%%%%%%%%%%%%%        20160922
{N}ote that $\psi(\theta_0)=1$. So
               for a triple $(p,q,\alpha)$
               if there exists some $\beta$ such that 
               \begin{equation}\label{2}
               \eta_\beta> 0 \ \mathrm{on\ } Q^\circ, % \mathrm{\ and\ } \lim_{\theta\rightarrow 0 \mathrm{\ or\ }\frac{\pi}{2}\psi\leq 1,
               \end{equation}
               then \eqref{1} holds on $Q$.
               \\{\ }
               
               \textbf{A1.} For simplicity, we try $\beta=1$ for {Row 1}. Now
               $$2\alpha={\Sigma}+2\ \mathrm{where}\ \Sigma=p+q,\ p=2r-2,\mathrm{\ and \ } q=2s-2.$$
                \begin{equation}\label{3}
               \begin{split}
                 4\alpha^2\eta_1  &=4\alpha^2-8(p+q+pq)+
                                        4(q^2-2q)\cot^2\theta+4(p^2-2p)\tan^2\theta\\
                                        &\geq 4\alpha^2-8(p+q+pq)+8\sqrt{(q^2-2q)(p^2-2p)}\ \ \ \ (\mathrm{``>"\ unless\ }p=q=2)\\
                                        &\geq 4\alpha^2-8(p+q+pq)+8(q-2)(p-2)\\
                                        &\geq (\Sigma+2)^2-24\Sigma+32\\
                                       % &=\Sigma^2-20\Sigma+36\\
                                        &=(\Sigma-2)(\Sigma-18)
                \end{split}
                 \end{equation}
Thus, when $2(r+s)-4=\Sigma\geq 18$, namely $r+s\geq 11$, we have $\eta_1>0$.
%%It can be checked that the second condition of \eqref{2} also holds for $r,s\geq 2$ with $r+s\geq 11$.
%Therefore \eqref{1} is true. 
%%$C_{r-1,s-1}$ for $r,s\geq 2$ with $r+s\geq 11$ are area-minimizing cones.

                         Actually from the first row of \eqref{3} one can get
                         \[
                         \begin{split}
                          4\alpha^2\eta_1 \tan^2\theta =&
                                    4(p^2-2p)\tan^4\theta+[4\alpha^2-8(p+q+pq)]\tan^2\theta+4(q^2-2q)\\
                                    =&4(p^2-2p)Y^2+[4\alpha^2-8(p+q+pq)]Y+4(q^2-2q)                                   
                          \end{split}
                          \]
                          where $Y=\tan^2\theta$. Easy to check that $\Delta$ of the quadratic is negative for $r,s\geq 3$ with $8\leq r+s\leq 10$.
                          So $\eta_1>0$. In particular $C_{2,4}$ for $(r,s)=(3,5)$
                          is area-minimizing. 
                          \\{\ }
                          
                          %%%%%%%
                          Due to Simons \cite{JS} and Simoes \cite{PS1,PS2} there are no area-minimizing cones for $r+s<8$ or for $(r,s)=(2,6)$.
                          Therefore $(r,s)=(2,7)$ and $(2,8)$ are the only left cases to check.
                          Although $df$ cannot be a global calibration on $Q$ for $(r,s)=(2,6)$, 
                          it is indeed a local calibration in a very {\it thin} angular neighborhood of the ray $\theta=\theta_0$, e.g. $1.145\leq\theta\leq 1.165$ shown in Figure 1 below.
                                     (At $\theta_0=\arctan \sqrt 5\approx 1.150262$, $\psi=1$.)
         In fact the existence of such angular neighborhood can be seen from $\eta_1(\theta_0)>0$.
                          As a consequence, $C_{1,5}$ is a stable minimal hypercone in $\mathbb R^8$.
          %%%
                \begin{figure}
                \begin{center}
                \includegraphics[scale=0.85]{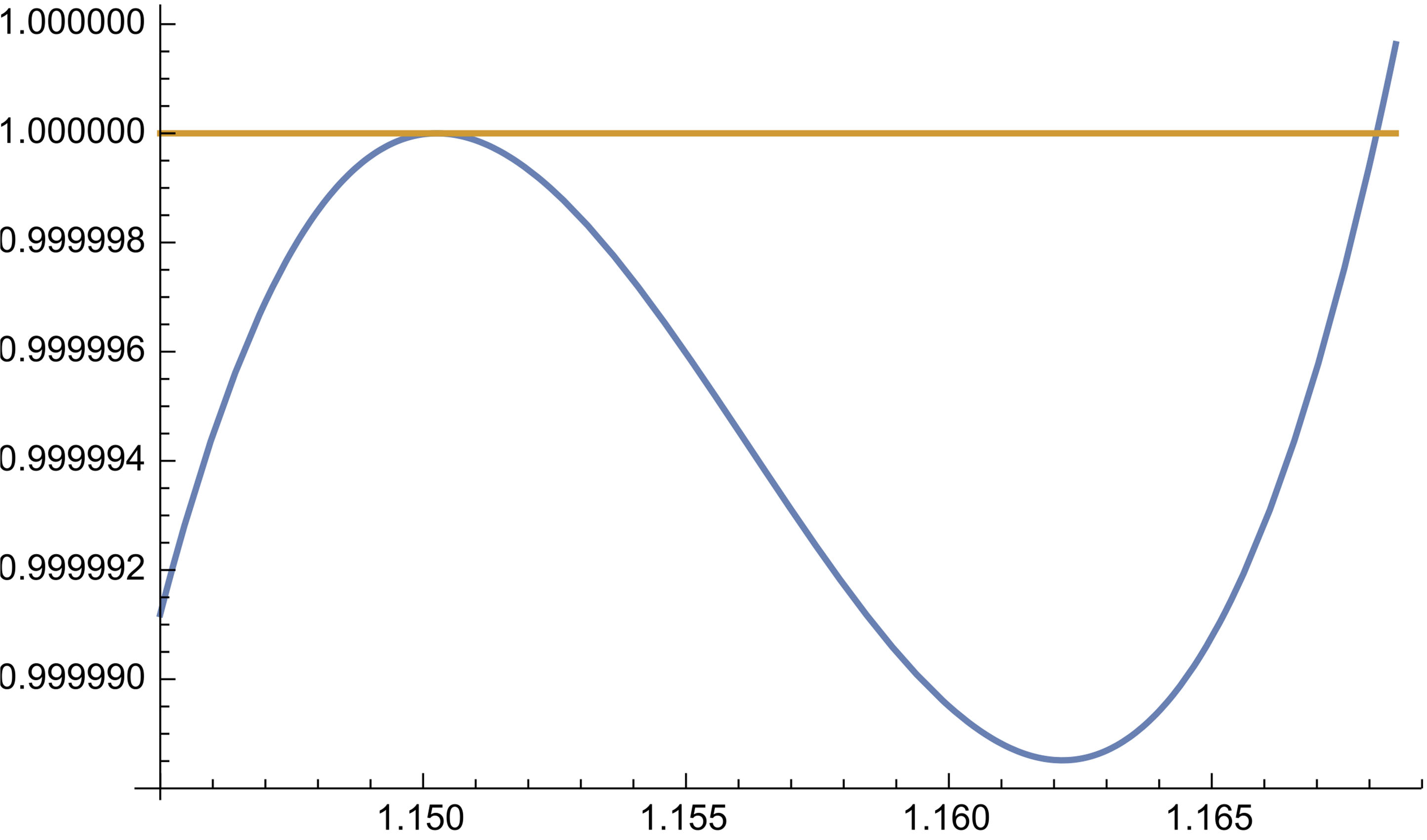}
                \caption{\small Graph of $\psi$ for $(r,s)=(2,6)$ with $\beta=1$.}
                \end{center}
                \end{figure}
%\\{\ }
  %   \\{\ }\\{\ }\\{\ }\\{\ }\\{\ }
                                      %%%%%%%
                           \textbf{A2.} For the rest two cases as well as for our purpose later in {C2.}, we consider $\beta=1.2>1$.
                        Similarly, it can be shown that $\eta_{1.2}>0$ when $r,s\geq 2$ and $r+s\geq 9$. 
                        \\{\ }
                                         
                        \textbf{B.} For 
                        $SO(2)\times SO(k)/\mathbb Z_2\times SO(k-2)$ in $\mathbb R^{2k}$ for $k=9$ (with $\alpha=8.5$)
                        one can use $\beta=1.2$.
                        There are two roots $\theta_{1,2}$ of $\eta_{1.2}$
                        satisfying $\theta_0<\theta_1<\theta_2$. 
                        Since $\psi(\theta_2)<1$, 
                        %we know that 
                        $df$ is a calibration on $Q$.
                        \begin{figure}
                        \begin{center}
                        \includegraphics[scale=0.85]{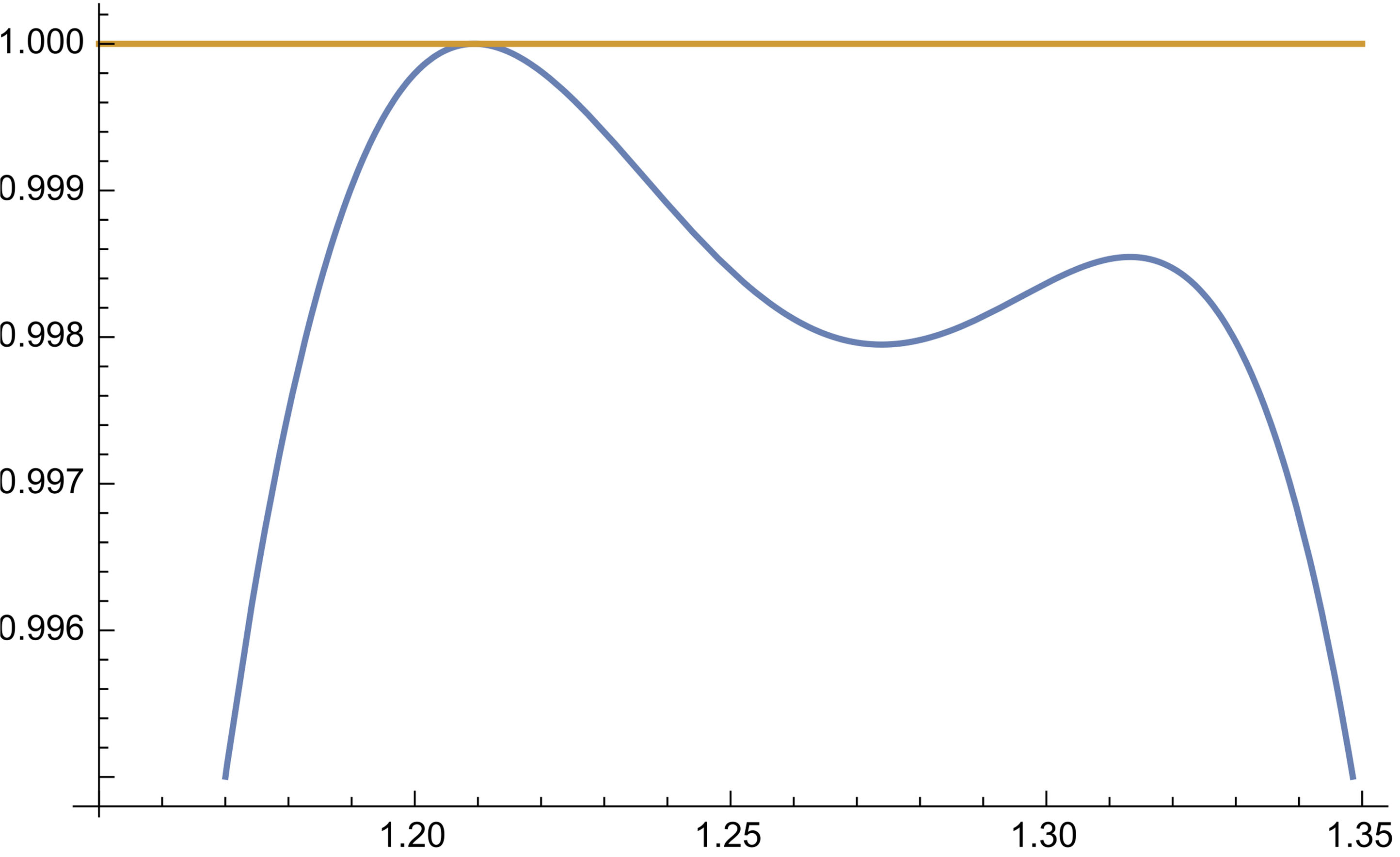}
                        \caption{\small Graph of $\psi$ with $\beta=1.2$.}                    
                        \end{center}
                        \end{figure}                    
%
                     %\\{\ }\\{\ }\\{\ }
                        
                        For $SO(2)\times SO(k)/\mathbb Z_2\times SO(k-2)$ in $\mathbb R^{2k}$ for $k=8$ (with $\alpha=7.5$),
                      using $\beta=1$,
                        there is only one root $\theta_1$ of $\eta_1$ and $\theta_1>\theta_0$.
                        So $df$ is a local calibration in some angular neighborhood of $\theta=\theta_0$, for example $|\theta-\theta_0|<\theta_1-\theta_0$.

                          \begin{figure}
                        \begin{center}
                        \includegraphics[scale=0.85]{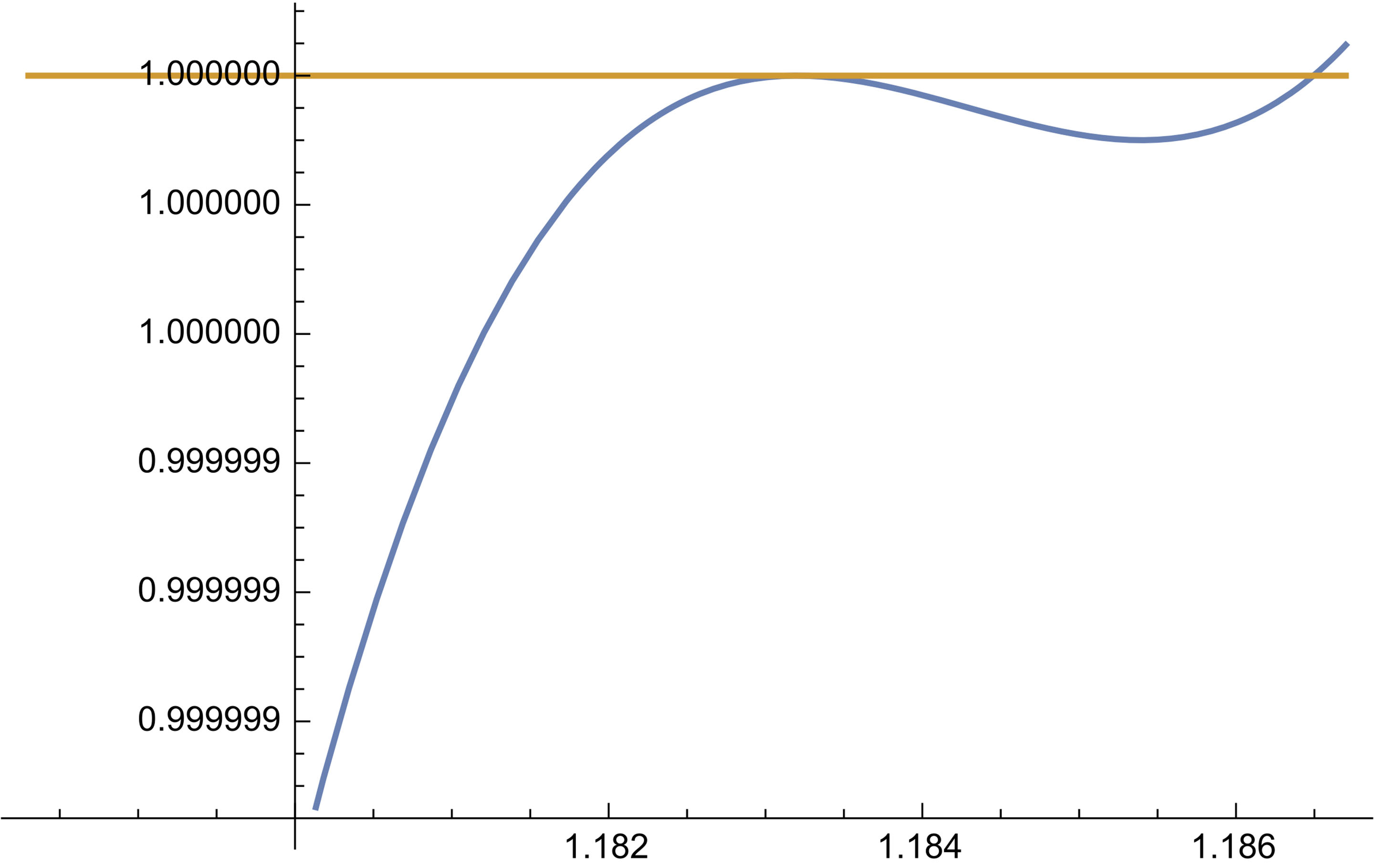}
                        \caption{\small Graph of $\psi$ with $\beta=1$.}                    
                        \end{center}
                        \end{figure}

                       Due to technique reason we save the the discussion on $SU(2)\times SU(k)/T^1\times SU(k-2)$ in $\mathbb R^{4k}$ for $k=4$ to \S\ref{s5}.

                        %%%%%%%%%%%%%%%%%%%%%%
                        {\ }
                        
                        \section{Proof of Theorem \ref{th:1.8}}\label{s4}
                        
                        By explanations on page 62 and Theorem 4.13 in \cite{HL2},
                                                 up to a constant,
                                                 \begin{equation}\label{CalinE}
                                                 \omega_f\triangleq \pi^*\left(df\right)\wedge \Omega,
                                                 \end{equation}
                                                 where $\Omega=\Omega_0/V$ and $\Omega_0$ is the oriented unit volume form of principal orbits,
                                                 is a coflat calibration of the cone, smooth away from singular orbits in the Euclidean space.
                                                 
                                                 In this section, we shall show that,
                                                 for an area-minimizing hypercone in  {Row 1},
                                                 either $\omega$
                                                 or its suitable modification is a calibration smooth away from the origin.  
                                                 
                                             Observe that the principal orbit through $(\overrightarrow X, \overrightarrow Y)\in\mathbb R^r\times\mathbb R^s$,
                                                 where neither $\overrightarrow X$ nor $\overrightarrow Y$ is a null vector,
                                                 is $S^{r-1}(\|\overrightarrow X\|)\times S^{s-1}(\|\overrightarrow Y\|)$.
                                                 Therefore 
                                                 \[
                                                 \Omega_0^{r,s}=\Omega_r(\overrightarrow X)\wedge\Omega_s(\overrightarrow Y), \mathrm{\ where}
                                                 \]
                                                 \[
                                                 \Omega_r(\overrightarrow X)=\frac{1}{\|\overrightarrow X\|}\sum_{i=1}^r  (-1)^{i-1}x_i dx_1\wedge\cdots\wedge \widehat{dx_i} \wedge\cdots\wedge x_r,
                                                \] 
                                                 \[
                     \Omega_s(\overrightarrow Y)=\frac{1}{\|\overrightarrow Y\|}\sum_{j=1}^s  (-1)^{j-1}y_j dy_1\wedge\cdots\wedge \widehat{dy_j} \wedge\cdots\wedge y_s,   
                                              \]%
                                              and
                                              $\{x_i\},\{y_j\}$ are standard basis of $\mathbb R^r$ and $\mathbb R^s$ respectively.
                        \\{\ }
                        
                        \textbf{C1.} When $r$ and $s$ are even integers with $r,s\geq 4$  in {A1.} ($\beta=1$),
                        there is a nonzero constant $C$ such that
                        \[
                        \begin{split}
                       &C\cdot\omega_f= C\cdot d\left(\pi^*f\right)\wedge\frac{\Omega_0^{r,s}}{V}\\
                =\ &d\left [\mathrm R^{3-r-s} \cdot \|\overrightarrow X\|^{2r-2} \cdot \|\overrightarrow Y\|^{2s-2} \right ]\wedge\frac{\Omega_0^{r,s}}{V}
                \ \ \ \left(\mathrm{\ Here\ R=\sqrt{\|\overrightarrow X\|^2+\|\overrightarrow Y\|^2}}.\ \right)\\
                =\ &
                      \Bigg \{\mathrm R^{3-r-s}\left[
                       \left(r-1\right)\|\overrightarrow X\|^{2r-4}\|\overrightarrow Y\|^{2s-2}d\|\overrightarrow X\|^2
                        +\left(s-1\right)\|\overrightarrow X\|^{2r-2}\|\overrightarrow Y\|^{2s-4}d\|\overrightarrow Y\|^2
                        \right]\\
                        &+\left(3-r-s\right)\mathrm R^{2-r-s}\|\overrightarrow X\|^{2r-2} \cdot \|\overrightarrow Y\|^{2s-2}d\mathrm R 
                      \Bigg \}\\
                    &\wedge \frac{1}{\|\overrightarrow X\|^{r}}\left(\sum_{i=1}^r  \left(-1\right)^{i-1}x_i dx_1\wedge\cdots\wedge \widehat{dx_i} \wedge\cdots\wedge x_r\right)\\
                   & \wedge\frac{1}{\|\overrightarrow Y\|^{s}}\left(\sum_{j=1}^s  \left(-1\right)^{j-1}y_j dy_1\wedge\cdots\wedge \widehat{dy_j} \wedge\cdots\wedge y_s\right)\\
                 =\ &\|\overrightarrow X\|^{r-4}\|\overrightarrow Y\|^{s-2}\cdot\Lambda_1+
                       \|\overrightarrow X\|^{r-2}\|\overrightarrow Y\|^{s-4}\cdot\Lambda_2+
                       \|\overrightarrow X\|^{r-2}\|\overrightarrow Y\|^{s-2}\cdot\Lambda_3\\
                        \end{split}
                        \]
                        where $\Lambda_i$ are\ smooth\ forms\ away\ from\ the\ origin.
                        Hence so is $\omega_f$.
                        \\{\ }
                        
                         \textbf{C2.} For remaining cases of {Row 1} in {A1.} and {A2.},
                         each reduced ray has a calibration $df$ where $f=C \alpha^{-1}r^{\alpha}c^{k}s^{t}$
                         for some positive constant $C$ and integers $k,t$ with $2k-p>2$ and $2t-q>2$.
                         (For such $f$ one can use $\beta=1$ when $p,q>2$ and $\beta=1.2$ when either $p$ or $q$ equals 2.)
                         Inspired by {C1.}, we wish to deform $f$ to $\tilde f$
                         so that $\omega_{\tilde f}$ behaves similarly as in {C1.} near singular orbits
                         and  $d\tilde f$ is a calibration of the cone.
                         
                         Our strategy is to increase the exponential of $f$.
                         We illustrate how to remove the singularity along $\theta=0$ (with $r>0$).                      
                         Choose an integer $N$ such that $N-\frac{q}{2}-1$
                         (which will replace the exponential $s-2$ of $\|\overrightarrow Y\|$ in {C1.})
                         is a positive even integer and $N>t$.
                         Then search for a suitable smooth function $ \lambda$ of $\theta$
                         for the deformation:
                         \begin{equation}\label{01}
                         \tilde f=C\alpha^{-1}r^{\alpha}c^ks^{t+\lambda}.
                         \end{equation}
                         We wish 
                               $\lambda(\theta)\equiv N-t$ when $\theta$ is very small
                               (for $\lambda=\lambda_{\epsilon_2}$ that we are about to construct,
                               such an interval of $\theta$ can be $(0,\frac{x_0}{2}-\epsilon_2)$ introduced below)
                               and
                               $ \lambda(\theta)\equiv0$ when $\theta>\tilde\theta$ for some small positive $\tilde \theta$.
                         By computation
                         \[
                         d\tilde f=C
{r^{\alpha-1}}  %%%%%%%%%%%%%%%%%%%%%%%%%%%%%%%%%%%%%%%%%%        20160922
                         c^ks^{t+\lambda}dr+C\alpha^{-1}r^{\alpha}c^ks^{t+\lambda}\left(\lambda'\ln s+\left(t+\lambda\right)\cot\theta
                         -k\tan\theta\right)d\theta
                         \]
                         and similar to \eqref{1} the square of comass of $d\tilde f$ equals
                         \begin{equation}\label{4}
                         \tau C^2c^{2k-p} s^{2t+2\lambda-q}\left[1+\left(\frac{1}{\alpha}\right)^2\left(\lambda'\ln s+\left(t+\lambda\right)\cot\theta-k\tan\theta\right)^2\right].
                         \end{equation}
                         In view of the square of comass of calibration $df$, 
                         \begin{equation*}\label{4.5}
  \tau C^2c^{2k-p} s^{2t-q}\left[1+\left(\frac{1}{\alpha}\right)^2\left(t\cot\theta-k\tan\theta\right)^2\right]\leq 1\mathrm{\ {on}\ } Q, %%%%%%%%%%%%%%%%%%%%%%%%%%%%%%%%%%%%%%%%%%        20160922
                         \end{equation*}
                         one only needs to show \eqref{4} is no large than one on the support of $\lambda$ for $d\tilde f$ being a calibration. 
                         
                         We consider the following {Lipschitzan} function for a fixed $0<x<\tilde \theta$:
                         \begin{equation*}\label{5}
                         %\displaystyle
                        {
                         \lambda_x(\theta)=
                                       \begin{cases}
                                     N-t                                   & 0\leq  \theta\leq \frac{x}{2}\\
                                     \left(N-t\right)
                                     \left[
                                     1-\dfrac{1}{
                                     {\int}^x_{\frac{x}{2}}\frac{1}{\ln \sin\gamma}\,d\gamma} \cdot {\int}_{\frac{x}{2}}^\theta\dfrac{1}{\ln \sin\gamma}\, d\gamma
                                     \right]
                                     \ \ \  \ \ \ \ \ & \frac{x}{2} <  \theta < x\\
                                     0                                     & x\leq\theta\\
\end{cases}}
                         \end{equation*}
                        By this choice \eqref{4} on $\frac{x}{2}<\theta<x$ becomes
                        \begin{equation}\label{6}
                         \tau C^2c^{2k-p} s^{2t+2\lambda_x-q}\left[1+\left(\frac{1}{\alpha}\right)^2
                         \left(-\dfrac{N-t}{{\int}^x_{\frac{x}{2}}\frac{1}{\ln \sin\gamma}\,d\gamma}+\left(t+\lambda_x\right)\cot\theta-k\tan\theta\right)^2\right].
                         \end{equation}
                         Focus on the behavior of \eqref{6} as $x$ tends to zero.
                         By the assumption $2t+2\lambda_x-q\geq 2\rho\triangleq 2t-q>2$,
                         the expression limits to (for $\tilde C^2=\frac{\tau C^2}{\alpha^2}$)
                       
                         \begin{equation}\label{lmt}
                        \lim_{x\rightarrow 0}\tilde C^2 s^{2t+2\lambda_x-q}
                         \left(\dfrac{N-t}{{\int}^x_{\frac{x}{2}}\frac{1}{\ln \sin\gamma}\,d\gamma}\right)^2\leq
                         \tilde C^2 (N-t)^2\left(\lim_{x\rightarrow 0}\dfrac{\sin^\rho x}{{\int}^x_{\frac{x}{2}}\frac{1}{\ln \sin\gamma}\,d\gamma}\right)^2.
                         \end{equation}
                         Since ${\int}_0^a\dfrac{1}{\ln \sin\gamma}\,d\gamma$ is convergent for every $0<a<\frac{\pi}{2}$,
                         the denominator of the latter limit in parentheses goes to zero and
                         one can apply  L'H$\mathrm{\hat{o}}$pital's rule.
                         Due to
                         $$\lim_{x\rightarrow 0}\dfrac{\rho (2\sin\frac{x}{2}\cos\frac{x}{2})^{\rho-1}\cos x}
                              {\frac{1}{\ln (2\sin\frac{x}{2}\cos\frac{x}{2})}-\frac{1}{\ln \sin\frac{x}{2}}\ \ }
                         =\lim_{x\rightarrow 0}
                         \dfrac{\rho (2\sin\frac{x}{2}\cos\frac{x}{2})^{\rho-1}\cos x}
                              {\frac{-\ln 2-\ln \cos\frac{x}{2}}{\ln (2\sin\frac{x}{2}\cos\frac{x}{2})\cdot\ln \sin\frac{x}{2}}}
                         =0    
                         $$
                         the left hand side of \eqref{lmt} is zero (strictly less than one).
                         
                         Not hard to see there exists a smooth approximation $\lambda$ of some $\lambda_x$ such that
                         \eqref{4} does not exceed one.
                         One construction can be given as follows.
                        Choose $x_0$ small enough so that, for some small positive $\epsilon_1<\frac{x_0}{5}$,
                         the following analogous expression to \eqref{6}
                         \begin{equation}\label{7}
                         \tau C^2c^{2k-p} s^{2\rho}\left[1+\left(\frac{1}{\alpha}\right)^2
                         \left(-\chi\cdot\dfrac{N-t}{{\int}^{x_0}_{\frac{x_0}{2}}\frac{1}{\ln \sin\gamma}\,d\gamma}+K\cot\theta-k\tan\theta\right)^2\right]
                         \end{equation}
                          is less than half of one pointwise for $\frac{x_0}{2}-\epsilon_1\leq\theta\leq x_0+\epsilon_1$, any $\chi\in[0,2]$ and $K\in[0,N+1]$. 
                         Let $\lambda_\epsilon$ be an integral convolution of $\lambda_{x_0}$
                                 by a compactly supported,  smooth, {\it even} mollifier function with averaging radius $\epsilon$ ($<\epsilon_1$).
                         Define $\sigma^\epsilon$ by
                         $$
                         \lambda_\epsilon'(y)=
                        - \sigma^\epsilon(y)
                        \dfrac{N-t}{{\int}^{x_0}_{\frac{x_0}{2}}\frac{1}{\ln \sin\gamma}\,d\gamma} \cdot \frac{1}{\ln \sin y} \mathrm{\ \ \ for\ }
                         \dfrac{x_0}{2}-\epsilon_1\leq y\leq x_0+\epsilon_1.
                        $$
                      Then $\sigma^\epsilon$ is a smooth function with 
                      \[
                      \lim_{\epsilon\rightarrow0}\sigma^\epsilon(y)=
                      \begin{cases}
                      1& \mathrm{\ for\ } \dfrac{x_0}{2}\leq y\leq x_0\\
                      \dfrac{1}{2} &\mathrm{\ at\ } x_0 \mathrm{\ or\ } \dfrac{x_0}{2}\\
                      0 & \mathrm{\ on\ } [\dfrac{x_0}{2}-\epsilon_1,\dfrac{x_0}{2})\bigcup(x_0,x_0+\epsilon_1]
                      \end{cases}
                      \]
                      Hence through a contradiction argument         
                                    there exists some small 
                     $\epsilon_2(<\epsilon_1)$
                     such that
                      $\sigma^{\epsilon_2}(y)\leq2 \mathrm{\ on\ } [\dfrac{x_0}{2}-\epsilon_1,x_0+\epsilon_1]$.
                      Comparing \eqref{4} and \eqref{7}
                      shows that $\lambda=\lambda_{\epsilon_2}$ meets our needs.
                    
                    According to the discussion in {C1.} and \eqref{CalinE}, $N-\frac{q}{2}-1$ being a positive even integer asserts the smoothness of $\omega_{\tilde f}$
                    along $\big\{(\overrightarrow X,\overrightarrow Y)\in\mathbb R^r\times\mathbb R^s:
                    \overrightarrow X=\overrightarrow 0 \mathrm{\ and\ }\overrightarrow Y\neq\overrightarrow 0\big\}$.
                    The same trick applies for removing the singularity
                    along $\big\{(\overrightarrow X,\overrightarrow Y)\in\mathbb R^r\times\mathbb R^s:
                    \overrightarrow X\neq\overrightarrow 0 \mathrm{\ and\ }\overrightarrow Y=\overrightarrow 0\big\}$.
                    Hence we finish the proof.
                    
                   {\ }
                    
                    \section{The Case of $SU(2)\times SU(k)/T^1\times SU(k-2)$ in $\mathbb R^{4k}$ for $k=4$}\label{s5}
                    In virtue of {\it Mathematica} it can be observed that the construction \eqref{0}
                    cannot directly produce a calibration for any value $\beta$ in this case. 
                    Instead we fix $\beta=1$ and look at the graph of 
                    $\psi=\frac{1}{\tau}
                    c^4s^{10} \left[1+\left(\frac{1}{7.5}\right)^2\left(10\cot\theta-4\tan\theta\right)^2\right]$.
                    %%%%%%%%%
                    \begin{figure}
                    \begin{center}
                    \includegraphics[scale=0.85]{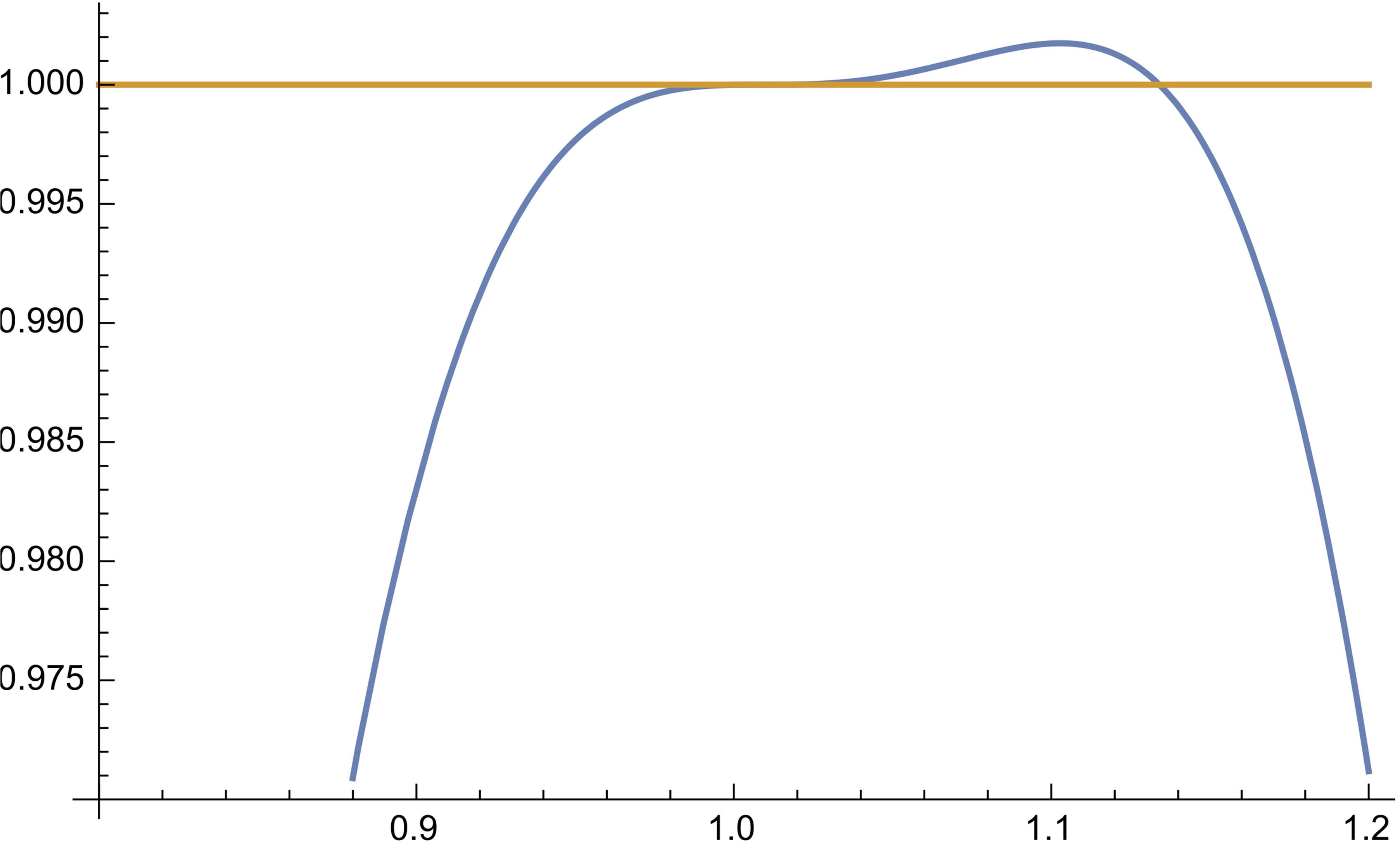}
                    \end{center}
                    \caption{}
                    \end{figure}
                                       Since the second peak is relatively not quite high, we try to chop it off.
                                       To do so, we alter the order of term $c$ by setting
                                                                 $$f_1=\frac{1}{7.5\tau}r^{7.5}c^{4+\lambda_1}s^{10},$$
                                     where $\lambda_1$ is a function. Then the square of the comass of $df_1$ becomes
                                       \begin{equation*}
                         (\star)\ \ \ \dfrac{1}{\tau}c^{4+2\lambda_1} s^{10}\left[1+\left(\frac{1}{7.5}\right)^2\left(\lambda_1'\ln c+10\cot\theta-(4+\lambda_1)\tan\theta\right)^2\right].
                         \end{equation*}  
                         Notice that $\theta_0$ (at where the first peak occurs) is less than 1.00686
                         and its local behavior is shown as follows:
                         
                         \begin{figure}%\label{F1}
                         \begin{center}
                    \includegraphics[scale=0.82]{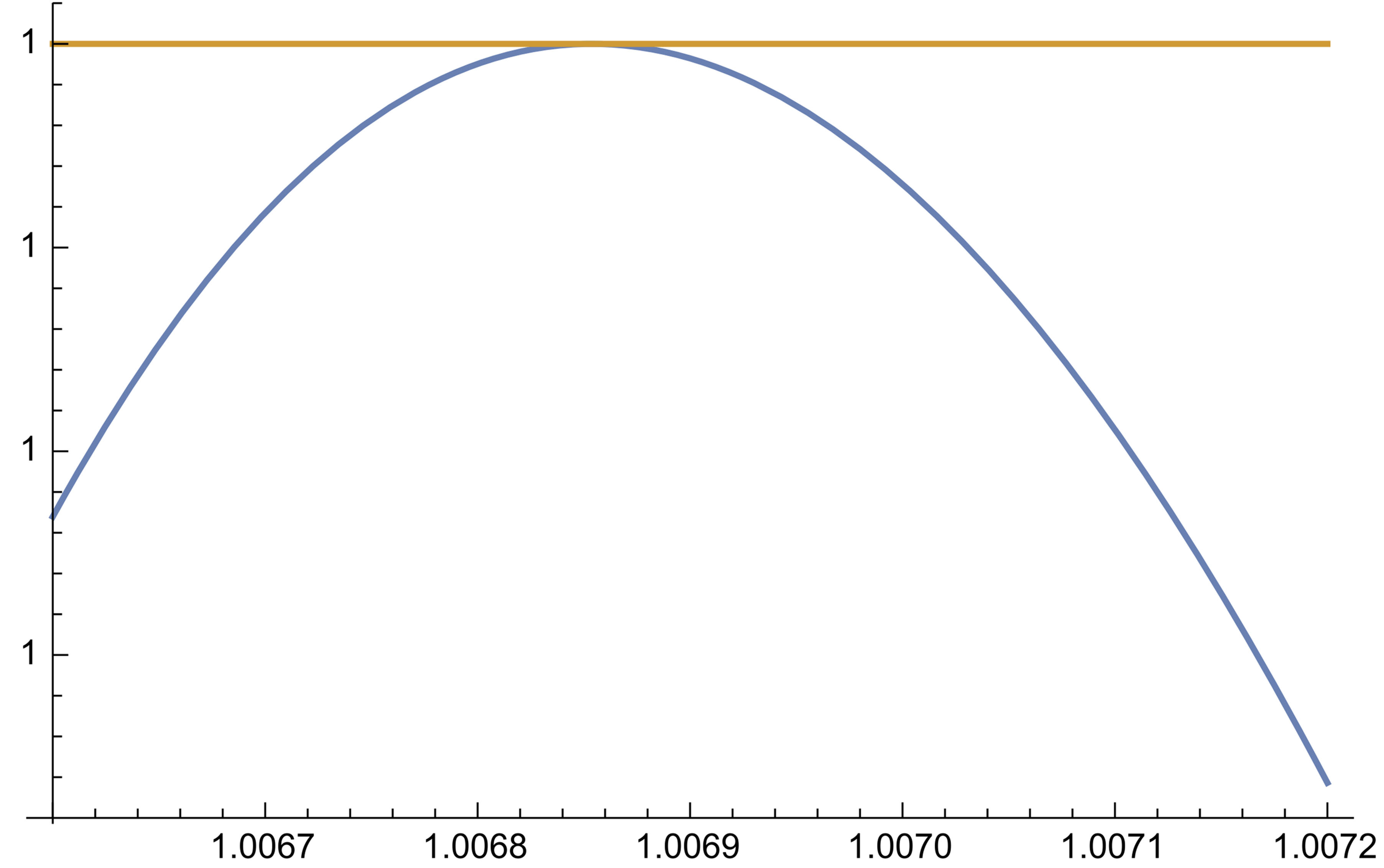}                   
                    \end{center}
                    \caption{}
                     \end{figure}

                         Now we consider the ODE: 
                         $$\lambda_1' = \dfrac{7.5\sqrt{\dfrac{1}{c^{4 + 2\lambda_1} s^{10} (\frac{14}{4})^2 (\frac{14}{10})^5} - 1}+ 
    10\cot\theta- \left(4 +\lambda_1\right) \tan\theta}{-\ln c},$$ with $\lambda_1(1.007) = 0.$
   Its solution tells us how to vary the exponential of $c$ to keep the corresponding comass exactly one on $[1.007,1.20]$.
   Running the following codes in {\it Mathematica}
                       {\small \[
                        \begin{split}
                                                \mathrm{sl =}&\mathrm{ NDSolve[\{y'[x] ==} \\
                                                &\mathrm{((7.5)(1/(Cos[x]^\wedge(4 + 2 y[x]) Sin[
                                                                     x]^\wedge(10) (14/4)^\wedge2 (14/10)^\wedge5) - 1)^\wedge(0.5) }\\
                                                &\mathrm{+ 
                                                                            10 Cot[x] - (4 + y[x]) Tan[x])/(-Log[Cos[x]]), y[1.007] == 0\}, 
                                                                              y, \{x, 1.007, 1.2\}]}
                                                                                                                  \end{split}
                                                 \]   
                            \ \ \ \ \ \ \    \ \ \ \     \         $\mathrm{Plot[Evaluate[y[x] /. sl], \{x, 1.007, 1.2\}]}$                                                  
                        }\\ 
                        we figure out the graph of $\lambda_1$.                                    
                        \begin{figure}
                                         
                                                                   \begin{center}
                    \includegraphics[scale=0.85]{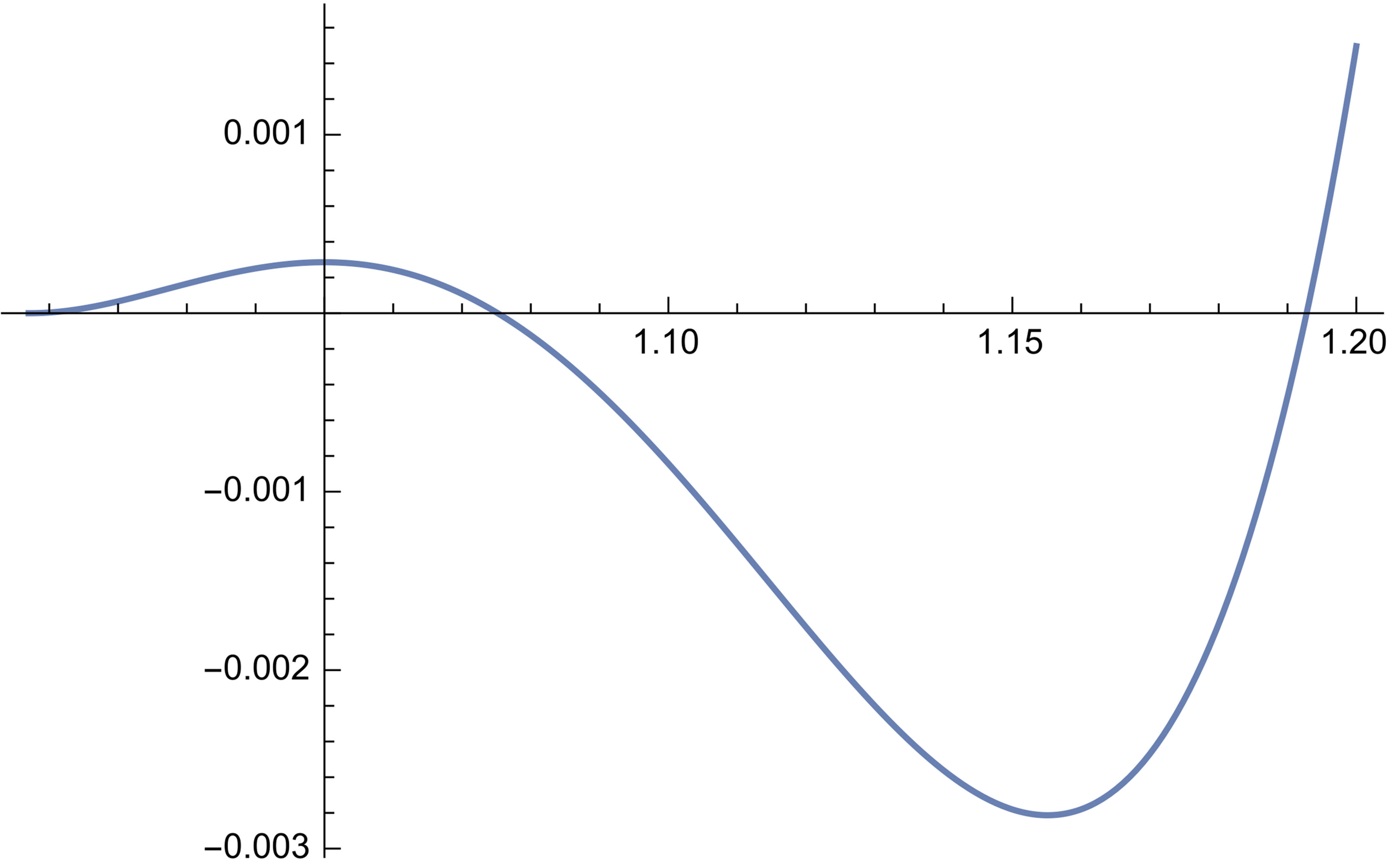}
                    \end{center}
                    \caption{}
                     \end{figure}
                    {\ } \\ 
                    Denote by $\theta_1$ the zero point of $\lambda_1$ near 1.2.
                    Note that  around 1.007
                    the integral curves $\Lambda$ of the slope field and the horizontal lines
                    form a coordinate chart.
                    Let the integral curve through $(1.007,0)$ be $\Lambda_{1.007}$.
                    Then one can glue 0 and $\Lambda_{1.007}$
                    such that on the left side of $\Lambda_{1.007}$
                    it transversally intersects each of $\Lambda$ at only one point.
                    Do the same kind of gluing around $\theta_1$.
                     \begin{figure}
                              \begin{minipage}[c]{0.5\textwidth}
                              \includegraphics[scale=0.75]{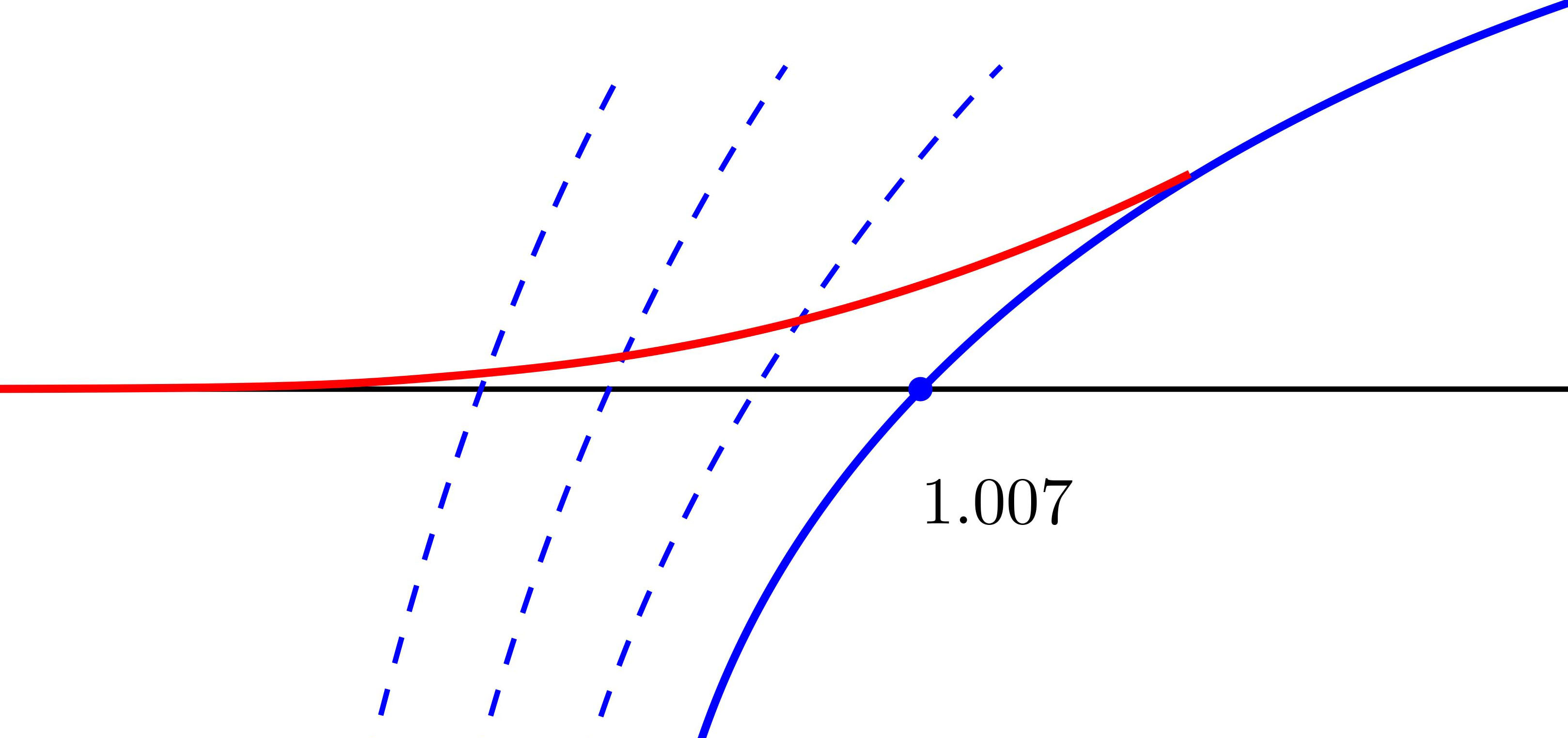}
                              \end{minipage}%
\ \ \ \  \
                          \begin{minipage}[c]{0.3\textwidth}
                           \includegraphics[scale=0.62]{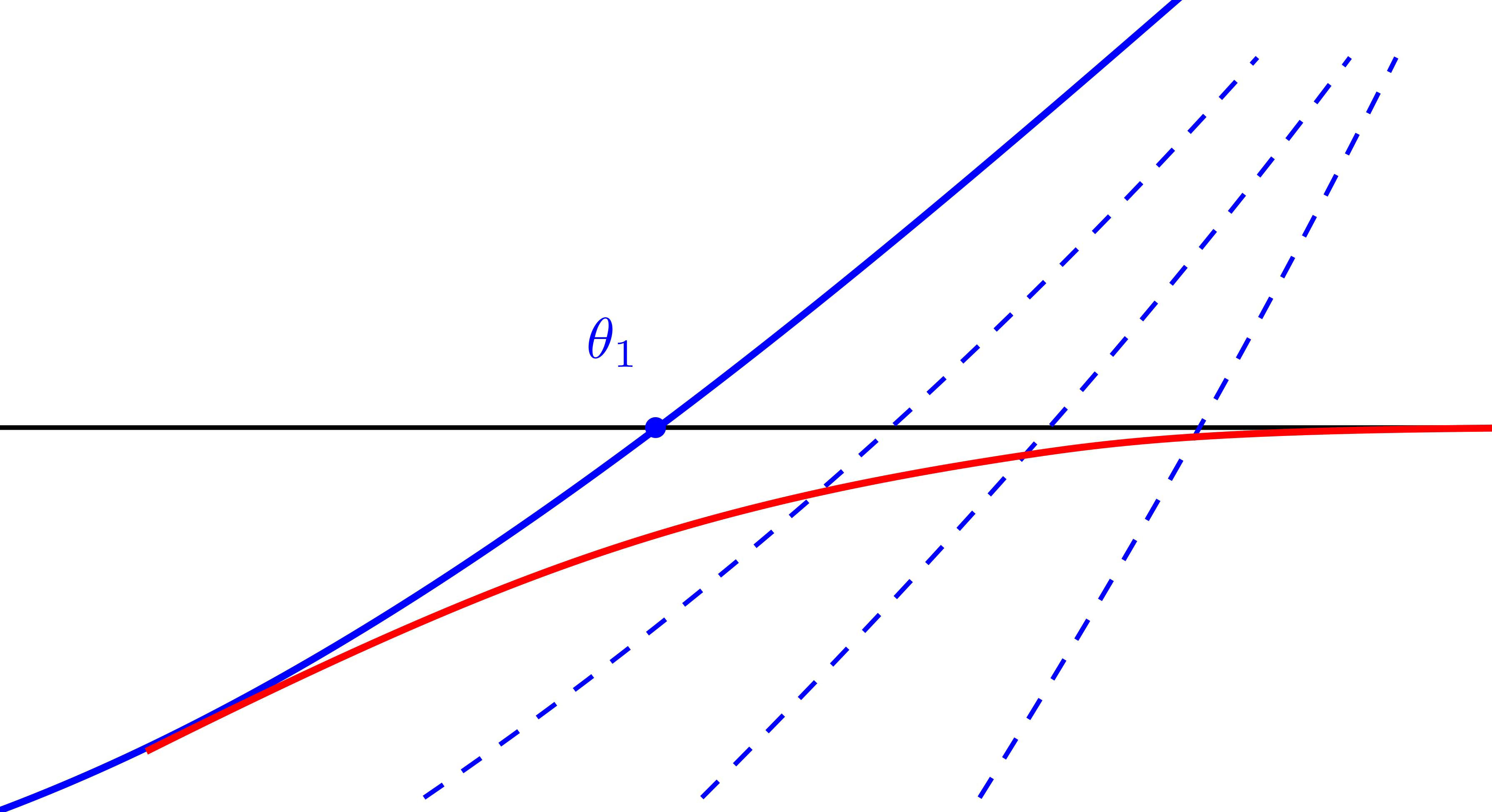}
                           \end{minipage}
                           \caption{}
                    \end{figure}
                    Now we get a smooth variation $\lambda$ of $\lambda_1$
                    that vanishes near $1.00686$ and $1.2$.
                    Since $10\cot\theta -4\tan\theta <0$ on $(\theta_0,1.2]$,
                                     after replacing $\lambda_1$ by $\lambda$, $(\star)$ is bounded from above by one.
                                     Hence we complete the proof of Theorem \ref{th:1.4}.

%{\ }\\%{\ }\\
     %%%%%%%%%%%%%%%%
                    \section{Proof of Theorem \ref{th:1.9}}
                    
                    Let us get back to the original calibration O.D.E. (or O.D.InE.)
                    on the (stretched) orbit space $Q$
                    of angle $\frac{\pi}{2}$
                    for area-minimizing hypercones of type \textbf{(II)}.
                    Set 
                    \[
                    f=\alpha^{-1}r^\alpha y(\theta)
                    \]
                    similar to \eqref{0}, then \eqref{code} becomes (w.r.t. $d\tilde s^2$)
                    \begin{equation}\label{Code}
                    y^2+(\frac{y'}{\alpha})^2\leq\frac{c^ps^q}{\tau}.
                    \end{equation}
                     \textbf{Remark.\ }
                    For Row 9-10 of type \textbf{(I)} with 
                         $V^2=c\ \cdot $ Im $\{(x+iy)^p\}^q$ and $A=\frac{\pi}{p}$.              
                      On the (stretched) orbit space $Q$
                    of angle ${\pi}$ with respect to (up to a constant factor)
                    \[
                    ds^2=r^{2\alpha-2}s^q(r^2d\theta+dr^2)\ \ \ \text{where}\ \ \ 
                    2\alpha-2=q-\frac{2p-2}{p},
                    \]
                       the corresponding calibration O.D.E. is
                    \begin{equation}\label{Code2}
                    y^2+(\frac{y'}{\alpha})^2\leq s^q.
                    \end{equation}
                    In this case $\theta_0=\frac{\pi}{2}$.
                                See \cite{BL} for details.
                    {\ }\\{\ }

                    Take type \textbf{(II)} for example (exactly the same for \textbf{(I)}).
                    For each area-minimizing hypercone of \textbf{(II)},
                    \eqref{Code} has a smooth solution $\Phi$
                    on $I=(0,\frac{\pi}{2})$ %for Row 1-6 %or on $I=(0,\pi)$ for Row 9-10
                    with value one at $\theta=\theta_0$ and value zero at the ending points of $I$.
                   Consider the slope fields $\Gamma_{1}, \Gamma_{2}$ of
                    \begin{equation}\label{y1}
                    y_1'=\alpha\sqrt{\frac{c^ps^q}{\tau}-y_1^2} \ \ \ \ \text{on}\ \ \ (0,\theta_0)
                    \end{equation}
                    and
                    \begin{equation}\label{y2}
                    y_2'=-\alpha\sqrt{\frac{c^ps^q}{\tau}-y_2^2} \ \ \ \ \text{on}\ \ \ I-(0,\theta_0].
                    \end{equation}
                    
                    We observe
                    (A) along the graph of $\Phi$, $\Gamma_1$ is positively deeper than $\Phi'$
                    and $\Gamma_2$ is negatively deeper than $\Phi'$
                    and
                    (B) for a fixed value $\theta<\theta_0$ (or $>\theta_0$),
                    the absolute value of the slope of $\Gamma_1$ (or $\Gamma_2$)
                    is decreasing in the variable $y$ (in the region below $y=\sqrt{\frac{c^ps^q}{\tau}}$).
                    
                     \begin{figure}
                                         
                                                                   \begin{center}
                    \includegraphics[scale=0.8]{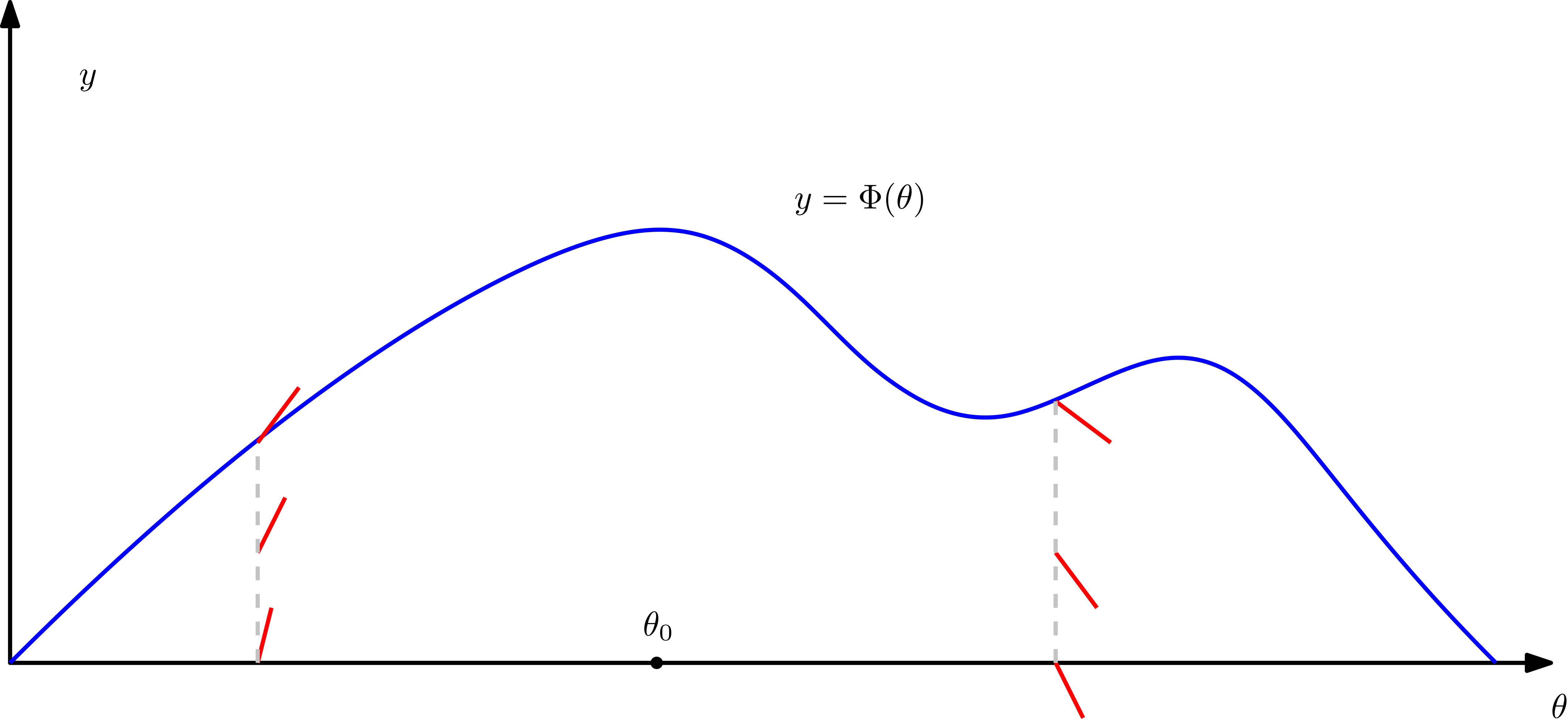}
                    \end{center}
                    \caption{}
                     \end{figure}

                    Therefore apparently there exists some nonnegative smooth solution $\Phi_0$ of \eqref{Code}
                    %as shown %in the picture
                    which is zero in some neighborhoods
                    around the ending points of $I$ and has value one at $\theta_0$.
                    
                     \begin{figure}
                                         
                                                                   \begin{center}
                    \includegraphics[scale=0.8]{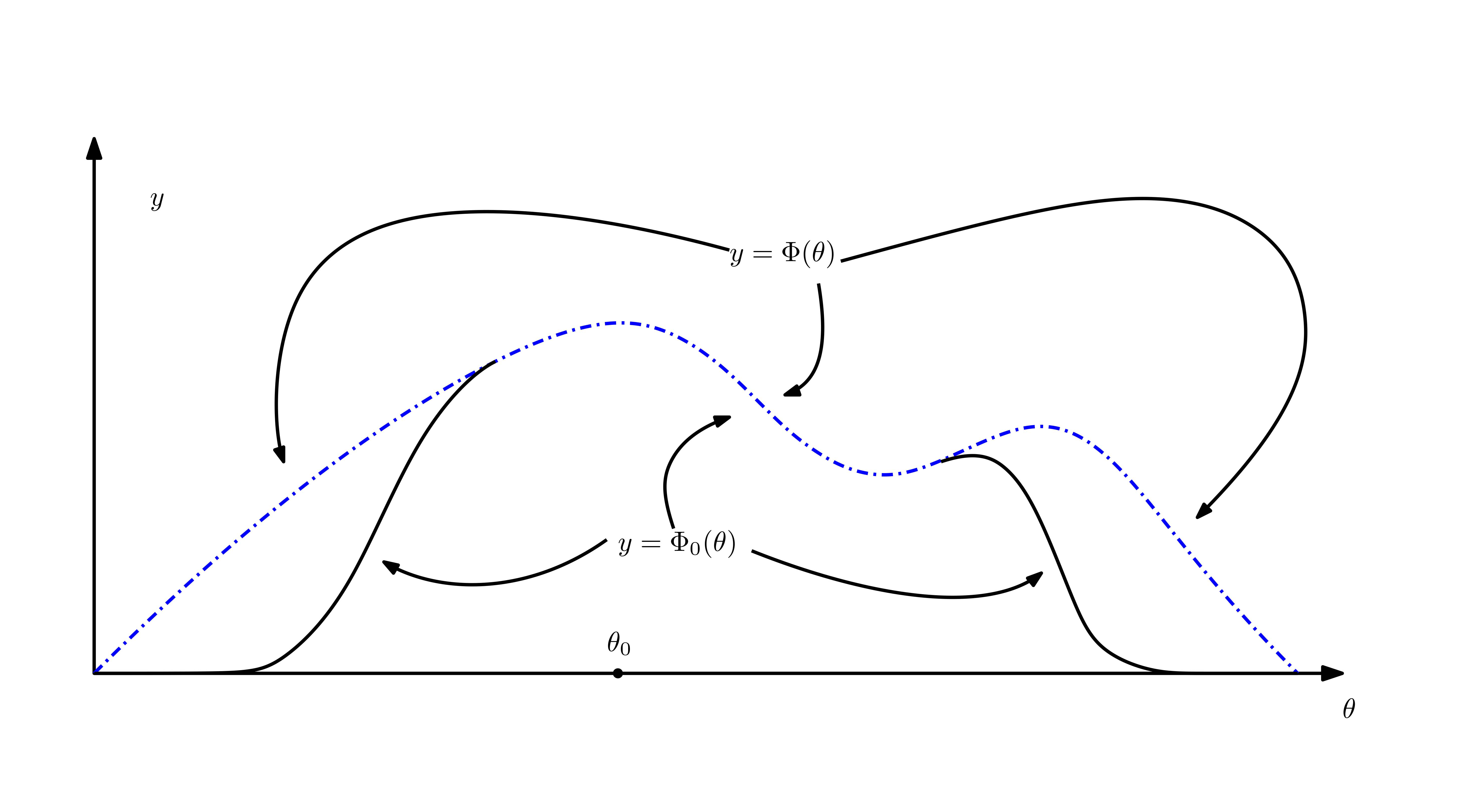}
                    \end{center}
                    \caption{}
                     \end{figure}

                    Then, according to \eqref{CalinE}, a coflat calibration singular only at the origin
                    can be gained for each homogeneous area-minimizing hypercone. 
                    Now the proof is complete.

{\ }

\acknowledgements{\rm The author would like to express deep gratitude to Professor H. Blaine Lawson, Jr. for his guidance and constant encouragement.
He also wishes to thank Professor Hui Ma and Doctor Chao Qian for helpful comments, the referees for nice suggestions,
and the MSRI for its warm hospitality.}

\end{document}